\documentclass{amsart}

\usepackage{amsmath,amsthm}%,amsfonts}
\usepackage {latexsym}
\usepackage{amssymb}

\newcommand{\less}{\lesssim}

\newcommand{\beal}{\begin{align}}
\newcommand{\enal}{\end{align}}
\newcommand{\bealn}{\begin{align*}}
\newcommand{\enaln}{\end{align*}}
\newcommand{\bear}{\begin{eqnarray}}
\newcommand{\eear}{\end{eqnarray}}
\newcommand{\beeq}{\begin{equation}}
\newcommand{\eneq}{\end{equation}}

\newcommand{\eps}{{\varepsilon}}
\newcommand{\R}{{\mathbb R}}

\newcommand{\calF}{{\mathcal F}}

\newcommand{\calS}{{\mathcal S}}

\newcommand{\calM}{{\mathcal M}}

\newcommand{\la}{\langle}
\newcommand{\ra}{\rangle}

\def\pr{\partial}

\def\nn{\nonumber}
\def\bm{\left[ \begin{array}{cc}}
\def\endm{\end{array}\right]}

\def\f{\frac}
\def\te{\tilde{e}}

\newtheorem{theorem}{Theorem}
\newtheorem{lemma}[theorem]{Lemma}

\newtheorem{cor}[theorem]{Corollary}
\newtheorem{prop}[theorem]{Proposition}

\theoremstyle{remark}

\def\il{\int\limits}

\renewcommand{\Re}{\,{\rm Re}}
\renewcommand{\Im}{\,{\rm Im}}
\def\slashint{-\!\!\!\!\!\!\il}

\begin{document}

\title {A remark on Littlewood-Paley theory for the distorted Fourier transform}
\author{W.\ Schlag}\thanks{The author was partially supported by the NSF grant
DMS-0300081 and a Sloan Fellowship. He wishes to thank Sigmund
Selberg from NTNU in Trondheim, Norway, for his hospitality and the
support by the Research Council of Norway through the research
project "Partial Differential Equations and Harmonic Analysis".
Also, he wishes to thank Fritz Gesztesy for many helpful
discussions.}

\address{Department of Mathematics, University of Chicago, 5734 South University Ave., Chicago, IL 60637, USA
} \email{schlag@math.uchicago.edu}

\maketitle

\section{Introduction}
\label{sec:intro}

The fundamental Littlewood-Paley theorem states that for any choice
of $1<p<\infty$ there is a constant $C_{p,d}$ such that
\begin{equation}
\label{eq:LP} C_{p,d}^{-1}\, \|f\|_p \less \| Sf\|_p \less C_{p,d}\,
\|f\|_p
\end{equation}
 for all Schwartz functions $f\in\calS(\R^d)$, see Stein~\cite{Stein1}, \cite{Stein2}.
  Here
\[ Sf:= \Big( \sum_{j=-\infty}^\infty |\Delta_j f|^2\Big)^{\f12}
\]
is the Littlewood-Paley square function which is associated with a
dyadic partition of unity
\[ \sum_{j=-\infty}^\infty \hat{\psi}(2^{-j}\xi) =1 \quad \forall\;
\xi\ne0 \] and $\Delta_j f = (\hat{\psi}(2^{-j}\cdot)
\hat{f})^{\vee}$. Closely related is the Mikhlin multiplier theorem,
which states that for all $1<p<\infty$
\begin{equation}
\label{eq:mikhlin_intro}
 \|(\mu \hat{f})^{\vee}\|_p \le C_{p,d}\, \|f\|_p
\end{equation}
provided $\mu\in C^{d+2}(\R^d\setminus\{0\})$ satisfies the
derivative bounds
\[ |\pr^\alpha \mu(\xi)|\le C_\alpha\, |\xi|^{-|\alpha|}
\]
for all $0\le |\alpha|\le d+2$, say (less is needed). The point here
is that the kernel $\hat{\mu}$ is a Calderon-Zygmund kernel. In
particular, it decays like $|x|^{-d}$ at infinity (and in general no
better, which means that $p=1$ or $p=\infty$ are forbidden). To pass
from the multiplier theorem to the Littlewood-Paley square function
bound one checks that the random function
\[ \mu(\xi) := \sum_j \pm \hat{\psi}(2^{-j}\xi)
\]
where $\pm$ are i.i.d.\ symmetric Bernoulli, are Mikhlin multipliers
uniformly in the choice of the signs $\pm$. The bound~\eqref{eq:LP}
is then obtained by combining~\eqref{eq:mikhlin_intro} with the
Khintchin inequality.

To purpose of this note is to investigate to what extent these same
bounds also hold relative to the distorted Fourier transform (we
will restrict ourselves to $d=3$). The latter here refers to a basis
of functions which diagonalize the Schr\"odinger operator
$H=-\Delta+V$ on its continuous spectral subspace. More precisely,
we need to assume that the real-valued potential $V$ decays
sufficiently rapidly to guarantee that this distorted Fourier basis
for $H$ exists. This basis, which we denote by
$\{e(x,\xi)\}_{\xi\in\R^d}$, consists of solutions of the
Lippman-Schwinger equation
\begin{equation}\label{eq:lipp}
e(\cdot,\xi) = e_\xi - R_V^{-}(\xi^2)V e_\xi, \qquad
e_\xi(x)=e^{ix\cdot\xi}
\end{equation}
See Agmon~\cite{Ag} for conditions that guarantee the solvability of
this equation as well as for the Fourier expansion and inversion
theorem in this setting\footnote{We are not after rough or slowly
decaying potentials here. In fact, we will consider only bounded
potentials that decay like $|x|^{-3-\eps}$}. Denote the distorted
Fourier transform by \begin{equation} \label{eq:FT} \calF_V f (\xi)
:= (2\pi)^{-\frac32}\il_{\R^3} \overline{e(x,\xi)} f(x)\, dx
\end{equation} as well as for any bounded $\mu=\mu(\xi)$ define a
multiplier operator by \begin{equation} \label{eq:Mmu} \calM_\mu f =
\calF_V^{-1} (\mu \calF_V f)
\end{equation}
Let $\sum_{j} \psi(2^{-j}\xi)=1$ be the usual Littlewood-Paley
decomposition of unity and  define a square function
\begin{equation}
\label{eq:SH} S_Hf(x) = \Big(\sum_j |\psi(2^{-j}\sqrt{H})
f|^2(x)\Big)^{\f12}
\end{equation}
on the continuous subspace $L^2_c(\R^3)$ relative to~$H$. It is of
course natural to develop analogues of the classical Mikhlin and
Littlewood-Paley theorems in this perturbed setting. By the
fundamental work of Yajima~\cite{Y} the wave operators are bounded
on $L^p(\R^3)$, $1\le p\le \infty$, provided $V$ decays faster than
a fifth power and $H$ has neither an eigenvalue nor a resonance at
zero energy. The latter here refers to the fact that there are no
solutions of the equation $Hf=0$, $f\in L^{2,-\sigma}(\R^3):=\la
x\ra^{\sigma} L^2(\R^3)$, $\sigma>\f12$. Equivalently, it means that
the perturbed resolvents $(-\Delta + V -z)^{-1}$ remain bounded as
operators from $L^{2,\sigma}\to L^{2,-\sigma}$ with $\sigma>1$ as
$z\to0$ in the upper half-plane $\Im z>0$. Since the wave operators
intertwine $-\Delta$ and $H$ restricted to its continuous subspace,
we conclude from Yajima's work that \eqref{eq:LP}
and~\eqref{eq:mikhlin_intro} remain valid for $S_H$ and $\calM_\mu$
provided we are on the continuous subspace~$L^2_c(\R^3)$.

In some applications, however, operators $H$ arise that do exhibit
resonances and/or eigenfunctions at zero energy. One such example,
which served as motivation for this work, was encountered by Krieger
and the author in~\cite{KS}: As discovered by Aubin, the functions
\[ \phi(x,a) = (3a)^{\f14}(1+a|x|^2)^{-\f12} \]
defined in $\R^3$ for each $a>0$  solve
\[ -\Delta \phi(\cdot,a) -\phi^5(\cdot,a) =0 \]
Hence, \beeq\label{eq:H_def} H_a :=-\Delta - 5\phi^4(\cdot,a)\eneq
satisfies
\[ H_a \, \pr_a \phi(\cdot,a) =0, \quad H_a \, \nabla \phi(\cdot,a) =0\]
Since $|\pr_a \phi(x,a) |\asymp |x|^{-1}$ as $|x|\to\infty$, it
follows that $\pr_a \phi \in L^{2,-\f12-\eps}\setminus L^2(\R^3)$.
This means that $H_a$ has both a resonance and an eigenvalue at zero
energy. In~\cite{KS}, the functions $\phi(\cdot,a)$ were considered
as static solutions of the wave equation $\Box \psi-\psi^5=0$. They
are linearly unstable, but they admit stable manifolds of
codimension one, at least for radial data, see~\cite{KS}. In order
to show this, the following linear bound was used (amongst others):
Fix $a>0$ and restrict $H_a$ to the radial subspace $L^2_{\rm
rad}(\R^3)$. Then there exists $c_0\ne0$ so that
\beeq\label{eq:sin_H} \f{\sin(t\sqrt{H_a})}{\sqrt{H_a}} =
c_0\,\pr_a(\cdot,a)\otimes \pr_a(\cdot,a) + \calS_a(t) \eneq where
$\calS_a(t)$ is dispersive in the following sense:
\[ \|\calS_a(t)f\|_\infty \less t^{-1} \|f\|_{W^{1,1}(\R^3)} \]
However, Strichartz estimates are still lacking for $\calS_a(t)$.
This is the main motivation for this note.

What we develop here is a Littlewood-Paley theorem for $H_a$ (or
more general operators with radial potentials) restricted to the
subspace of radial functions. It turns out that due to the
singularity of the resolvent at zero energy which results from the
resonance, the Mikhlin multiplier theorem can only hold in the range
$\f32<p<3$. In this range, we show that it does hold at least for
radial functions. The reason for this smaller range of $p$ lies with
the decay of the kernel associated with $\calM_\mu$. It decays at
infinity like $|x|^{-2}$ rather than $|x|^{-3}$ as in the
Calderon-Zygmund case. This forces $p>\f32$ and by duality also
$p<3$. Hence, we can only prove the Littlewood-Paley theorem in this
same range of~$p$ (and for radial functions). This is inadequate for
deriving Strichartz estimates for $\calS_a(t)$. Consequently, it
seems reasonable to believe that the latter requires a
Littlewood-Paley theory which is based on the regular part of the
spectral measure of $H$, just as much as $\calS_a(t)$ is the sine
evolution of the regular part of the spectral measure. However, we
do not pursue that here.

A very detailed analysis of Besov spaces relative to the distorted
Fourier transform was carried out by Jensen and Nakamura~\cite{JN1},
\cite{JN2}. However, these authors only consider inhomogeneous Besov
spaces. In other words, they carry out the dyadic partition  only
for large energies and they treat small energies as a single block.
This is not only quite different from the full square function, but
is also insufficient for proving  Strichartz estimates for the
operator~$\calS_a(t)$.

\section{Littlewood-Paley theory in the perturbed radial case}
\label{sec:LP}

Let $H=-\Delta+V$ with a real-valued, radial, decaying potential $V$
in $\R^3$ and with no eigenvalue, but a resonance at
zero\footnote{Hence, the results from this section apply to $H_a$ as
in~\eqref{eq:H_def} provided we restrict $H_a$ to $L^2_{\rm
rad}(\R^3)$}. As for the decay, we shall assume that $|V(r)|\less
\la r\ra^{-3-\eps}$ with $\eps>0$. It will also be convenient to
assume that $|V'(r)|\less \la r\ra^{-4-\eps}$ with $\eps>0$.  In
this section, we shall only consider radial multipliers operating on
radial functions. Let $\{e(x,\xi)\}$ be the distorted Fourier basis
for $H$, see above. Let, for any $r,k>0$ and $x$ with $|x|=r$,
\begin{equation}\label{eq:etilde} \te(r,k) = \f{rk}{4\pi}\il_{S^2}
e(x,k\omega)\, \sigma(d\omega)\end{equation} which is well-defined
by the radial symmetry of $V$. Then, for any $k>0$,
\begin{equation}\label{eq:1dSch}
(-\pr_{rr} + V)\te (r,k) = k^2 \te(r,k) \quad \forall
r>0, \qquad \te(0,k)=0
\end{equation}
and, moreover, for any radial Schwartz functions $f,g$
\begin{align}
&\la \calM_\mu f,g\ra
= \la \calF_V^{-1} (\mu \calF_V f), g \ra = \nn \\
&=(2\pi)^{-3}\il_0^\infty\il_0^\infty\il_0^\infty
\mu(k)\il_{S^2}\il_{S^2}\il_{S^2} e(r_1\omega_1, k\omega)\,
\overline{e(r_2\omega_2,k\omega)}
\,\sigma(d\omega_1)\sigma(d\omega_2)\,\sigma(d\omega)\;
k^2\,dk\; r_1^2 f(r_1) r_2^2 g(r_2)\, dr_1 dr_2 \nn \\
&= 8\il_0^\infty\il_0^\infty\il_0^\infty \mu(k) \te(r_1, k)\,
\overline{\te(r_2,k)} \, dk\; r_1f(r_1) r_2 g(r_2)\, dr_1 dr_2 \nn \\
&=\f{2}{\pi}\il_0^\infty\il_0^\infty\il_0^\infty \mu(k) \te(r_1,
k)\, \overline{\te(r_2,k)} \, dk\; \tilde f(r_1) \tilde g(r_2)\,
dr_1 dr_2 \nn \\
&= \la \tilde\calF_V^{-1} (\mu \tilde\calF_V \tilde f), \tilde g \ra
=: \la \widetilde\calM_\mu \tilde f,\tilde g\ra \label{eq:radial}
\end{align}
where $\tilde f=\sqrt{4\pi}\,rf(r)$, $\tilde g=\sqrt{4\pi}\,rg(r)$
and\footnote{For a particularly simple derivation of the Fourier
expansion theorem in the half-line case see Gesztesy and
Zinchenko~\cite{GZ}.} \beeq\label{eq:FT1d} (\tilde \calF_V \tilde
f)(r) = \sqrt{\f{2}{\pi} }\il_0^\infty
\overline{\te(r,\rho)}\,\tilde f(\rho)\, d\rho\eneq
 This shows that it suffices to study spectral multipliers
associated with the operator on $L^2(0,\infty)$ given by
\[ \tilde H := -\pr_{rr} + V(r) \]
with a Dirichlet condition at $r=0$. Our normalizations are chosen
so that $f\mapsto \tilde f$ as a map $L^2_{\rm rad}(\R^3)\to
L^2(0,\infty)$ is unitary and so that both $\calF_V$ and $\tilde
\calF_V$ are unitary on their respective $L^2_c$ spaces. In the free
case of course $e(x,\xi)=e^{ix\cdot\xi}$ and $\tilde
e(r,k)=\sin(rk)$.

Our first goal in this section is to prove the following result.
Notice that we do not make any assumption on eigenvalues or
resonances at zero energy.

\begin{prop}\label{prop:1dmult}
Let $ H=-\Delta+V(r)$ with $|V(r)|\less \la r\ra^{-\beta}$ and
$|V'(r)|\less \la r\ra^{-\beta-1}$ for some $\beta>3$.  Let $\mu$ be
an arbitrary radial function in $C^3(\R^3\setminus\{0\})$
satisfying\footnote{Less is needed on $\mu$, but we ignore these
refinements} $|\mu^{(\ell)}(k)|\le
 k^{-\ell}$ for all $k>0$ and $\ell=0,1,2,3$. Then
$\calF_V^{-1} (\mu \calF_V f)=:\calM_\mu f$ satisfies the bounds
\beeq\label{eq:mikhlin} \|\calM_\mu f\|_p \le C(V,p) \|f\|_p \eneq
for all $3/2<p<3$ and all radial $f\in L^p(\R^3)$.
\end{prop}

The proof of this proposition exploits the reduction to the
half-line operator $\widetilde \calM_\mu$ from~\eqref{eq:radial}. It
will be based on the theory of $A_p$ weights on the line, see
Stein~\cite{Stein2}, Chapter~V. The main estimates will be given in
two separate pieces for large and small $k$, see
Lemma~\ref{lem:high} and~\ref{lem:low} below, respectively. We need
to develop some standard background material before embarking on the
proofs.

\noindent The solution to~\eqref{eq:1dSch} is unique up to
normalization. We remark that in our application~\eqref{eq:1dSch}
has a bounded, nonzero solution for $k=0$ due to the resonance at
zero. For arbitrary $k>0$,
\[ \te (r,k) = c_+(k) f(r,k) + c_-(k) \overline{f(r,k)} \]
where $f(r,k)$ are the Jost functions, i.e.,
\[ (-\pr_{rr} + V)f (r,k) = k^2 f(r,k) \quad \forall
r>0, \qquad f(r,k)\sim e^{irk} \text{\ \ as\ \ }r\to\infty
\]
They are the unique solutions of the Volterra equation
\[ f(r,k) = e^{irk} + \il_r^\infty \f{\sin(k(r'-r))}{k} V(r')
f(r',k)\, dr' \qquad \forall r\ge0 \] If $k=0$, then
\[ f(r,0) = 1 + \il_r^\infty (r'-r) V(r')
f(r',0)\, dr' \qquad \forall r\ge0 \]
 It is standard that
\[ \sup_{r,k\ge 0}|f(r,k)|\le \exp\Big( \int_0^\infty r'|V(r')|\, dr'\Big)<\infty \]
The resolvent of $-\pr_{rr}$ on $L^2(0,\infty)$ with Dirichlet
boundary condition is
\[ R_0^+(k^2)(r,r') = k^{-1} \sin(rk) e^{ir'k} \chi_{[0<r<r']} +
k^{-1} \sin(r'k) e^{irk} \chi_{[0<r'<r]}
\]
Therefore, by the Lippman-Schwinger equation for $\tilde H$,
\begin{align*} \te(r,k) &= \sin(rk) - R_0^+(k^2) V
\te(\cdot,k)(r) \\
&= \sin(rk) - k^{-1} \il_0^r \sin(r'k) e^{irk} V(r') \te(r',k)\, dr'
-k^{-1} \il_r^\infty \sin(rk) e^{ir'k} V(r') \te(r',k)\,dr'
\end{align*}
which implies that
\begin{align*}
\pr_r \te(0,k) &= k - \il_0^\infty e^{ir' k} V(r') \te(r',k)\, dr'
\\
& = k - c_+(k) \il_0^\infty e^{ir'k} V(r') f(r',k)\, dr' - c_-(k)
 \il_0^\infty e^{ir'k} V(r') \overline{f(r',k)}\, dr'
 \end{align*}
Hence $c_+(k), c_{-}(k)$ solve the system
\begin{align*}
c_+(k)f(0,k) + c_-(k) \overline{f(0,k)} & = 0 \\
c_+(k)\Big(\pr_r f(0,k)+ \il_0^\infty e^{ir'k} V(r') f(r',k)\,
dr'\Big) + c_-(k) \Big(\overline{\pr_r f(0,k)}+\il_0^\infty e^{ir'k}
V(r') \overline{f(r',k)}\, dr'\Big) & = k
\end{align*}
Let $D(k)$ denote the determinant of this system. Then
\[
D(k) = f(0,k)\overline{\pr_r f(0,k)} - \overline{f(0,k)}\pr_r f(0,k)
+ \il_0^\infty e^{ir'k} V(r') [f(0,k)\overline{f(r',k)} -
\overline{f(0,k)} f(r',k)]\, dr'
\]
Since
\begin{align} f(0,k) &= 1 + k^{-1}\il_0^\infty \sin(r'k)
V(r') f(r',k)\, dr' \label{eq:f0}\\ f'(0,k) = \pr_r f(0,k) &= ik -
\il_0^\infty \cos(r'k) V(r') f(r',k)\, dr' \nn\end{align} we further
conclude that
\begin{align*} D(k) &= -2ik\Re (f(0,k)) + i \il_0^\infty \sin(r'k)
V(r') 2i\Im
[f(0,k)\overline{f(r',k)}]\, dr' \\
&=-2ik\Re (f(0,k)) + 2 k\Im [f(0,k)(1-\overline{f(0,k)})] = -2ik
{f(0,k)}
\end{align*}
Evaluating the Wronskian $W[f(\cdot,k),\overline{f(\cdot,k)} ]$ at
$r=0$ and $r=\infty$ implies that
\begin{equation}\label{eq:wronsk}
\Im(f(0,k)\overline{f'(0,k)}) = -k
\end{equation} and thus
\beeq\label{eq:f0k} C(V)
> |f(0,k)| > C(V)^{-1} k(1+k)^{-1}
\eneq
 Finally, for all $k\ne0$
\[ c_+(k) = \f{1}{2i} \f{\overline{f(0,k)}}{f(0,k)}, \qquad c_-(k) =
-\f{1}{2i}
\]
and thus
\begin{equation}\label{eq:te_form} \te(r,k) =
\f{f(r,k)\overline{f(0,k)}- \overline{f(r,k)}f(0,k)}{2if(0,k)}
\end{equation}
Since $f(r,k)=e^{irk} + O(k^{-1})$, it follows that
\[ \te(r,k)=\sin(rk) + O(k^{-1}) \text{\ \ as\ \ }k\to\infty \]
The behavior at $k=0$ depends on whether or not there is a nonzero
bounded solution $y$ to $\tilde H y=0$, $y(0)=0$. It is easy to see
that this dichotomy is equivalent with the dichotomy $f(0,0)=0$
vs.~$f(0,0)\ne0$. Indeed, since $|V(x)|\less \la x\ra^{-3-\eps}$,
all such solutions $y$ are given by
\[ y(r) = c_1 + c_2r + \il_r^\infty (r'-r) V(r') y(r')\, dr' \]
Hence, $y\in L^\infty(0,\infty)$ iff $c_2=0$ in which case
$y(r)=c_1f(r,0)$. Hence, there exists a nonzero bounded solution $y$
to $\tilde H y=0$, $y(0)=0$ iff $f(0,0)=0$. In our
application~\eqref{eq:H_def}, $H_a(\pr_a \phi(\cdot,a))=0$ and thus
$\tilde H_a(r\pr_a \phi(\cdot,a))=0$ with a bounded solution
$y(r)=r\pr_a(r,a)$ vanishing at $r=0$. So there $f(0,0)=0$, which
implies that~\eqref{eq:f0k} is optimal in that case.

In what follows, we will assume that $f(0,0)=0$. This is the harder
case, and $f(0,0)\ne0$ is also implicit in what we are doing below.
Hence, differentiating~\eqref{eq:wronsk} in~$k$ then proves that
\[ \Im (\pr_k f(0,0)\overline{f'(0,0)}) = -1 , \qquad \pr_k f(0,0)\ne0\]
and by d'Hospital's rule $c_+(k)$ is continuous on $[0,\infty)$ with
\[ c_+(0)= \f{1}{2i} \f{\overline{\pr_k f(0,0)}}{\pr_k f(0,0)} \]
This requires that $\la r\ra^2 V\in L^1$, which is the case here.
However, we cannot guarantee that $c_+\in C^1[0,\infty)$ since that
would require $\la r\ra^3 V\in L^1$. Nevertheless, in view of
\eqref{eq:f0} one has $f(0,k)\in C^\infty(0,\infty)$ and thus also
$c_+\in C^\infty(0,\infty)$.

The following lemma reduces $\widetilde \calM_\mu$ to manageable
pieces in the high-energy case.

\begin{lemma}
\label{lem:high} Let $\mu\in C^2(0,\infty)$ with
$|\mu^{(\ell)}(k)|\le  k^{-\ell}$ for all $k>0$ and $\ell=0,1,2$.
Assume further that $\mu(k)=0$ for all $0<k<1$.  Then the kernel
$K(r,r')$ of  $\widetilde\calM_\mu=\tilde\calF_V^{-1} \mu
\tilde\calF_V $ satisfies
\[ K = K_1+K_2+K_3 \]
where \begin{align*} K_1(r,r') &= \f12 \int_0^\infty \cos((r-r')k)\,
\mu(k)\, dk \\
|K_2(r,r')| &\less \la r-r'\ra^{-2} \\
|K_3(r,r')| &\less (r+r')^{-1}
\end{align*}
for all $r,r'\in (0,\infty)$. In particular, all these kernels are
$L^p(0,\infty)$-bounded for $1<p<\infty$.
\end{lemma}
\begin{proof}
It will be convenient to assume that $\mu(k)=0$ for large values of
$k$, but the bounds will not depend on this additional cut-off so
that it can be removed in the end. The kernel of $\calM_\mu$ is
\begin{align}
K(r,r') &= \il_0^\infty \mu(k) \te(r,k) \overline{\te(r',k)}\, dk \nn\\
&= \il_0^\infty \mu(k) \Big(c_+(k)f(r,k)+c_-(k)\overline{f(r,k)}\Big)
\Big(\overline{c_+(k)}\,\overline{f(r',k)}+\overline{c_-(k)}f(r',k)\Big)\, dk \nn \\
&= \f14\il_0^\infty e^{i(r-r')k}\mu(k) m(r,k)\; \overline{m(r',k)}
\, dk+\f14\il_0^\infty e^{-i(r-r')k}\mu(k) m(r',k)\;
\overline{m(r,k)} \,
dk \label{eq:r-r'}\\
&\quad +\f1{2i}\il_0^\infty e^{i(r+r')k}\mu(k) c_+(k)
m(r,k){m(r',k)} \, dk-\f1{2i}\il_0^\infty e^{-i(r+r')k}\mu(k)
\overline{c_+(k)} \;\overline{m(r,k)\,m(r',k)} \, dk \label{eq:r+r'}\\
&=: K^{(+,+)}(r,r') + K^{(-,-)}(r,r') +K^{(+,-)}(r,r') +
K^{(-,+)}(r,r')\nn
\end{align}
where we have set $f(r,k)=e^{irk} m(r,k)$. Thus,
\beeq\label{eq:m_form} m(r,k) = 1 + \il_r^\infty
\f{e^{2ik(r'-r)}-1}{2ik}\, V(r') m(r',k)\, dr'\eneq which shows that
\[ \sup_{r>0}|m(r,k)-1|\less k^{-1} , \qquad \sup_{r,k>0} |\pr_r
m(r,k)|<\infty\]
 Using a derivative of $V$, we can improve on the first bound. Indeed,
\[ m(r,k) = 1 - \il_r^\infty \f{e^{2ik(r'-r)}-1}{(2ik)^2}\, \big[V'(r')
m(r',k)+ V(r')\pr_{r'} m(r',k)\big]\, dr' \] and thus $m_1:=m-1$
satisfies \beeq\label{eq:m1_est}  \sup_{r>0}|\pr_k^j m_1(r,k)|\less
k^{-2} \qquad j=0,1,2 \eneq Therefore, also
\[ \sup_{r,r'>0}|\pr_k^j (m(r,k)m(r',k)-1)|\less k^{-2} \qquad j=0,1,2 \]
which implies that
\[ \Big|K^{(+,+)}(r,r') - \f14\il_0^\infty e^{i(r-r')k}\mu(k)\,
dk\Big|\less \la r-r'\ra^{-2} \] Since it is standard that
\[\il_0^\infty e^{i(r-r')k}\mu(k)\, dk \] is a singular integral
kernel, $K^{(+,+)}$ satisfies the desired $L^p$-bounds, and so does
$K^{(-,-)}$. To analyze $K^{(+,-)}$ and $K^{(-,+)}$, we first note
that
\[ c_+(k) = \f{1}{2i} \f{1+\overline{m_1(0,k)}}{1+ m_1(0,k)} \]
so that by \eqref{eq:m1_est}
\[ \big|\pr_k^j\;[ \mu(k) c_+(k)
m(r,k){m(r',k)}]\big |\less k^{-j} \qquad j=0,1,2 \] Hence,
\[ |K^{(\pm,\mp)}(r,r')|\less (r+r')^{-1} \]
Since this latter kernel is well-known to be $L^p(0,\infty)$ bounded
for $1<p<\infty$, we are done.
\end{proof}

Next, we consider small energies.

\begin{lemma}
\label{lem:low} Let $\mu\in C^3(0,\infty)$ with
$|\mu^{(\ell)}(k)|\le k^{-\ell}$ for all $k>0$ and $\ell=0,1,2,3$.
Assume further that $\mu(k)=0$ for all $k>1$. Then the kernel
$K(r,r')$ of  $\widetilde\calM_\mu=\tilde\calF_V^{-1} \mu
\tilde\calF_V $ satisfies\footnote{If $f(0,0)=0$, then $K_1=0$}
\[ K = K_1+K_2+K_3 \]
where \begin{align*} K_1(r,r') &= \f12 m(r,0)m(r',0)\int_0^\infty
\cos((r-r')k)\,
\mu(k)\, dk \\
|K_2(r,r')| &\less |r-r'|^{-1},\quad |\pr_r K_2(r,r')|+|\pr_{r'} K_2(r,r')|\less |r-r'|^{-2} \\
|K_3(r,r')| &\less (r+r')^{-1}
\end{align*}
for all $r,r'\in (0,\infty)$. All these kernels are bounded on
$L^p(0,\infty)$ for $1<p<\infty$.
\end{lemma}
\begin{proof}
We will again rely on the decomposition \eqref{eq:r-r'}
and~\eqref{eq:r+r'}. With $\psi$ our Littlewood-Paley function, we
further write
\begin{align}\label{eq:Ksplit} K^{(+,+)}(r,r') &=
\f14\il_0^\infty e^{i(r-r')k} \mu(k) m(r,0)\;\overline{m(r',0)} \,
dk +\sum_{j<0} \f14\il_0^\infty e^{i(r-r')k} \psi(2^{-j}k)\mu(k)
m(r,r';k)
 \, dk \\
 &=: K^{(+,+)}_0(r,r') + \sum_{j<0} K^{(+,+)}_j(r,r') \nn
 \end{align}
with
\[ m(r,r';k) = m(r,k)\;\overline{m(r',k)}-m(r,0)\;\overline{m(r',0)}=m(r,k)\;\overline{m(r',k)}-m(r,0)\;{m(r',0)}\]
Since $m(\cdot,k)\in L^\infty (0,\infty)$, the first integral
$K^{(+,+)}_0$ on the right-hand side of~\eqref{eq:Ksplit} is bounded
$L^p(0,\infty)\to L^p(0,\infty)$ for $1<p<\infty$. Here we are again
using that
\[ \il_0^\infty
e^{i(r-r')k} \mu(k)\, dk\] is a singular integral kernel. To obtain
$K_1(r,r')$, define
\[ K_1(r,r') = K^{(+,+)}_0(r,r') + K^{(-,-)}_0(r,r') \]
 We now
claim that \beeq\label{eq:claim} \Big|\pr_r K^{(+,+)}_j(r,r')\Big|
\less \min(2^{2j}, |r-r'|^{-3} 2^{-j}) \qquad \forall j<0,\;\forall
r,r'>0\eneq and symmetrically with $\pr_{r'}$. This in turn implies
that $\tilde K^{(+,+)}(r,r'):=\sum_{j<0}K^{(+,+)}_j(r,r')$ satisfies
\[ \Big|\pr_r \tilde
K^{(+,+)}(r,r')\Big| \less |r-r'|^{-2} \] from which we obtain the
H\"ormander condition \beeq\label{eq:hor}
\sup_{r_1,r_2>0}\il_{[|r'-r_1|>2|r_1-r_2|]}| \tilde
K^{(+,+)}(r_1,r') - \tilde K^{(+,+)}(r_2,r')|\, dr' < \infty \eneq
To prove \eqref{eq:claim}, let $h(a)=\f{e^{2ia}-1}{2ia}$. This
function satisfies $|h^{(\ell)}(a)|\less \la a\ra^{-\ell-1}$.
Then~\eqref{eq:m_form} is the same as
\[ m(r,k) = 1 + \il_r^\infty h(k(r'-r)) (r'-r) V(r') m(r',k)\, dr'
\]
and thus \beeq\label{eq:prkm}
 \pr_k m(r,k) = \il_r^\infty h'(k(r'-r)) (r'-r)^2 V(r') m(r',k)\,
 dr' + \il_r^\infty h(k(r'-r)) (r'-r) V(r') \pr_k m(r',k)\, dr'
 \eneq
 as well as
\begin{align} \pr_r\pr_k m(r,k) &= -k\il_r^\infty h''(k(r'-r)) (r'-r)^2 V(r')
m(r',k)\,
 dr' - 2\il_r^\infty h'(k(r'-r)) (r'-r) V(r')  m(r',k)\, dr'\nn\\
& - \il_r^\infty e^{2ik(r'-r)} V(r') \pr_k m(r',k)\,
dr'\label{eq:prrkm}
\end{align}
From these identities we obtain the bounds \beeq\label{eq:prm}
\sup_{0<k<1}\sup_{r>0}[|\pr_k\, m(r,k)|+|\pr_r\pr_k\,\, m(r,k)|]<
\infty \eneq and therefore also\footnote{Automatically in this proof
$0<k<1$}
\[ \sup_{r>0}\Big(|\pr_r(m(r,k)-m(r,0))|+ |m(r,0)-m(r,k)|\Big) \less k \]
which finally implies
\[ \sup_{r,r'>0}\Big(|m(r,r';k)|+|\pr_r m(r,r';k)|\Big) \less k \]
These bounds give us one half of the claim~\eqref{eq:claim}. Indeed,
\[
\Big|\pr_r K^{(+,+)}_j(r,r')\Big| \less \il_0^\infty |\psi(2^{-j}k)|
k\, dk \less 2^{2j}
\]
To obtain the other half, we estimate
\[
|r-r'|^3\big|\pr_r K^{(+,+)}_j(r,r')\big| \less \il_0^\infty \Big|
\pr_k^3[k\,\psi(2^{-j}k) \mu(k)\,m(r,r';k)]\Big|\, dk + \il_0^\infty
\Big| \pr_k^3[\psi(2^{-j}k) \mu(k)\,\pr_r\,m(r,r';k)]\Big|\, dk
\]
These estimates imply the remaining half of~\eqref{eq:claim}, viz.
\[|r-r'|^3\big|\pr_r K^{(+,+)}_j(r,r')\big| \less 2^{-j}\] provided we
can prove that\footnote{We only need this for $\ell\le3$, but it
holds for all $\ell$} \beeq \label{eq:mrr'k_bds}
\sup_{r,r'>0}[|\pr_k^\ell\, m(r,r';k)| + |\pr_k^\ell\pr_r
\,m(r,r';k)| ]\less k^{1-\ell} \qquad \forall \;\ell\ge0\eneq
Writing \begin{align*} m(r,r';k) &=
(m(r,k)-m(r,0))m(r',k)+m(r,0)(m(r',k)-m(r',0)) \\
\pr_r m(r,r';k) &=
(\pr_r\,m(r,k)-\pr_r\,m(r,0))m(r',k)+\pr_r\,m(r,0)(m(r',k)-m(r',0))
\end{align*}
it is easy to see that \eqref{eq:mrr'k_bds} follows from
\eqref{eq:prkm} and~\eqref{eq:prrkm}. For example,
\begin{align*}
& \pr^2_k m(r,k) = k^{-1}\il_r^\infty h''(k(r'-r))k(r'-r)\; (r'-r)^2
V(r') m(r',k)\,
 dr' \\
 &+ 2\il_r^\infty h'(k(r'-r)) (r'-r)^2 V(r') \pr_k m(r',k)\, dr'
   + \il_r^\infty h(k(r'-r)) (r'-r) V(r') \pr^2_k m(r',k)\, dr'
\end{align*}
proves that $|\pr^2_k m(r,k)|\less k^{-1}$ and inductively,
$|\pr^\ell_k m(r,k)|\less k^{1-\ell}$ for all $\ell\ge2$. Similarly,
\begin{align*}
\pr_r\pr_k^2\, m(r,k) &= -\il_r^\infty h''(k(r'-r)) (r'-r)^2 V(r')
m(r',k)\,
 dr' -\il_r^\infty h'''(k(r'-r)) k(r'-r)\;(r'-r)^2 V(r')
m(r',k)\,
 dr' \\
 & -k\il_r^\infty h''(k(r'-r)) (r'-r)^2 V(r')
\pr_k\,m(r',k)\,
 dr' - 2\il_r^\infty h''(k(r'-r)) (r'-r)^2 V(r') \pr_k m(r',k)\, dr'
 \\
 & - 2\il_r^\infty h'(k(r'-r)) (r'-r) V(r') \pr^2_k\, m(r',k)\, dr'
 - 2i\il_r^\infty e^{2ik(r'-r)} (r'-r) V(r') \pr_k m(r',k)\, dr' \\
 & + V(r) \pr_k m(r,k)
\end{align*}
implies that $|\pr_r\pr_k^2\, m(r,k)|\less k^{-1}$, and analogously
for the higher derivatives.

To summarize, \[ K_2(r,r')=\tilde K^{(+,+)}(r,r')+\tilde
K^{(-,-)}(r,r') \] satisfies
\[ |\pr_r K_2(r,r')|+ |\pr_{r'} K_2(r,r')|\less |r-r'|^{-2} \]
and therefore a H\"ormander condition. It is easier and also
implicit in the preceding to check that
\[ |K_2(r,r')\less |r-r'|^{-1} \]
We also need to estimate the kernel $K_3(r,r'):=K^{(+,-)}(r,r') +
K^{(-,+)}(r,r')$ from~\eqref{eq:r+r'}. For those we remark that the
bounds
\[ |\pr_k^\ell [\mu(k)c_+(k) m(r,k)m(r',k)]|\less k^{-\ell}\qquad
\ell=0,1,2 \] which follow easily from our preceding work, imply
\[ | K^{(+,-)}(r,r')| +|
K^{(-,+)}(r,r')| \less (r+r')^{-1}
\]
However, this latter kernel $(r+r')^{-1}$ is $L^p(0,\infty)$
bounded.

We yet need to establish the $L^p$ boundedness of $\tilde
K^{(+,+)}(r,r')+\tilde K^{(-,-)}(r,r')$. Now suppose that this sum
is an $L^2$ bounded operator. Since we have a H\"ormander condition,
the Calderon-Zygmund theorem (see Stein~\cite{Stein2}, page~19,
Theorem~3) will then guarantee that it is bounded on $L^p$ with
$1<p\le2$ and by taking adjoints, this will then also hold for $2\le
p<\infty$. Hence it suffices to check the case $p=2$. First,
$\tilde\calM_\mu$ with kernel $K(r,r')$ is $L^2$ bounded by the
spectral theorem. Second,
\[\tilde
K^{(+,+)}(r,r')+\tilde K^{(-,-)}(r,r') = K(r,r') - K^{(+,+)}_0(r,r')
- K^{(-,-)}_0(r,r') - K^{(+,-)}(r,r') - K^{(-,+)}(r,r') \] is
therefore also $L^2$ bounded. The lemma is proved.
\end{proof}

These two lemmas establish $L^p(0,\infty)$ boundedness of
$\widetilde\calM_\mu$ for any Mikhlin multiplier $\mu$ on the
half-line. However, we are really after the $L^p(\R^3)$-boundedness
of $\calM_\mu$. In view of~\eqref{eq:radial}, this is equivalent
with a certain weighted $L^p$ boundedness of $\widetilde\calM_\mu$.
The theory of $A_p$ weights will give us what we need, as we shall
see now. It is in order to apply the theory of $A_p$ weights that we
have expressed the kernels explicitly as in the previous two lemmas.

\begin{proof}[Proof of Proposition~\ref{prop:1dmult}]
First, observe that for radial functions $f,g$ in $\R^3$,
\[ \la \calM_\mu f,g \ra = \la \widetilde \calM_\mu \tilde f,\tilde
g\ra
\] and thus
\begin{align*}
\|\calM_\mu f\|_p &= \sup_{\|g\|_{p'}=1} \big|\la \widetilde\calM_\mu \tilde f,\tilde g\ra \big|\\
&= \sup_{\|r^{\f{2}{p'}-1} \tilde g\|_{L^{p'}(0,\infty)}=1}
\big|\la r^{1-\f{2}{p'}}\widetilde\calM_\mu \tilde f,r^{\f{2}{p'}-1}\tilde g\ra \big| \\
&= \big\| r^{\f{2}{p}-1}\widetilde\calM_\mu \tilde f \big
\|_{L^p(0,\infty)}
\end{align*}
Hence, we need to prove that
\[
\| r^{\f{2}{p}-1}\widetilde\calM_\mu \tilde f\|_p \less
\|r^{\f{2}{p}-1} \tilde f\|_p
\]
on the half-line. We now claim that $\omega(r)=r^{2-p}$ is an
$A_p$-weight for all $3/2< p<3$. Recall that this means that
\beeq\label{eq:Ap} \sup_{b>a\ge0}\quad \slashint_{a}^b \omega(x)\,
dx \;\Big( \slashint_a^b \omega^{-\f{p'}{p}}(x)\,
dx\Big)^{\f{p}{p'}} < \infty \eneq where $\slashint$ is the averaged
integral.  To check this for our weight $\omega(r)=r^{2-p}$, we
first remark that due to $3/2<p<3$, we have $1>(2-p)p'/p$ and
$2-p>-1$. Hence, if $b\gg a$, then the left-hand side
of~\eqref{eq:Ap} is
\[ \sup_{b>0} \less b^{2-p} \big(b^{-(2-p) p'/p}\big)^{p/p'} \less 1 \]
If $b\less a$, then one applies the mean-value theorem to reach a
similar conclusion.

The main property of $A_p$ weights is that the Hardy-Littlewood
maximal operator, as well as Calderon-Zygmund operators are bounded
on $L^p(\omega)$. Since one can check this manually for
$\omega(r)=r^{2-p}$ and the kernel $(r+r')^{-1}$, the proposition
now follows from Lemmas~\ref{lem:high} and~\ref{lem:low}.
\end{proof}

It is now standard to pass from Proposition~\ref{prop:1dmult} to the
Littlewood-Paley theorem. Let $\sum_{j} \psi(2^{-j}k)=1, \; k\ne0$
be the usual Littlewood-Paley decomposition of unity.

\begin{cor}
\label{cor:LP} Let $ H$ be as in Proposition~\ref{prop:1dmult} and
define
\[ S_{ H}f = \Big(\sum_j \Big|\psi\Big(2^{-j}\sqrt{ H}\Big) f\Big|^2\Big)^{\f12}
\]
to be the square function relative to $H$. Then
\begin{equation}\label{eq:LP_upbd}
 C(V,p)^{-1}\|f\|_p\le \| S_{H}f \|_p \le C(V,p)\|f\|_p \eneq for all
radial $f\in L^2_c(\R^3)\cap L^p(\R^3)$ with $3/2<p<3$, where
$L^2_c(\R^3)$ is the continuous subspace of $ H$.
\end{cor}
\begin{proof}
The upper bound follows from \eqref{eq:mikhlin} by means of the
usual randomization and Khinchin inequality method, whereas the
lower bound in~\eqref{eq:LP_upbd} then follows from the upper bound
via duality. We skip the details.
\end{proof}

\section{Concluding remarks}

\begin{itemize}

\item
It would be interesting to extend the results here to the case of
non-radial functions, as well as non-radial potentials. In this case
it is most likely necessary to estimate the solutions of the
Lippman-Schwinger equation directly. Here we chose to work in the
radial setting, since the ensuing reduction to the half-line allows
us to work with a Volterra integral equation, rather than the more
difficult  Fredholm integral equation (as in Lippman-Schwinger). The
Volterra equation of course is the one governing the Jost solutions.
At any rate, in the general case it will  be necessary to isolate
the singularity of the perturbed resolvent at energy zero as in
Jensen and Kato~\cite{JenKat}, or Jensen and Nenciu~\cite{JenNen}.
See also~\cite{ES} where the Jensen-Nenciu method is carried out
in~$\R^3$.

\item As mentioned before, Strichartz estimates for the wave equation with a potential
remain open for potentials that lead to singular resolvents at zero
energy (we are referring here to Strichartz estimates on the regular
part of the evolution, and wish to obtain the full range of
exponents as in the free case). The same appears to be true for the
Schr\"odinger case. However, that case should be much easier since
it does not require the Littlewood-Paley theorem. For the case of
resonance/eigenvalue at zero energy the {\em dispersive} estimates
for the Schr\"odinger equation with a potential are in~\cite{ES} as
well as Yajima~\cite{Yaj2}, whereas for the wave equation with a
potential they are derived in~\cite{KS} (although in the latter case
only a resonance is allowed). The latter paper motivated the
question of Strichartz estimates for the dispersive part of the
evolution operator (i.e., $\calS_a(t)$ in~\eqref{eq:sin_H}).

\end{itemize}

\bibliographystyle{amsplain}

\begin{thebibliography}{99}

\bibitem[Agm]{Ag} Agmon, S. {\em Spectral properties of Schr\"odinger
operators and scattering theory.} Ann.\ Scuola Norm.\ Sup.\ Pisa
Cl.\ Sci. (4) 2 (1975), no.~2, 151--218.

\bibitem[ErdSch]{ES} Erdo\smash{\u{g}}an, M. B., Schlag, W.
{\em Dispersive estimates for Schr\"{o}dinger operators in the
presence of a resonance and/or an eigenvalue at zero energy in
dimension three: I}, Dynamics of PDE, vol.~1, no.~4 (2004),
359--379.


%\bibitem[Gol]{Gold} Goldberg, M. {\em Dispersive bounds for the
%three-dimensional Schr\"odinger equation with almost critical
%potential}, preprint 2004, to appear in GAFA.

\bibitem[GesZin]{GZ} Gesztesy, F., Zinchenko, M. {\em On spectral
theory of Schr\"odinger operators with strongly singular
potentials.} preprint 2005.

\bibitem[JenKat]{JenKat} Jensen, A., Kato, T. {\em Spectral properties of
Schr\"odinger operators and time-decay of the wave functions.} Duke
Math.\ J.\ 46  (1979), no. 3, 583--611.


\bibitem[JenNak1]{JN1} Jensen, A., Nakamura, S. {\em $L^p$ and Besov estimates for
Schr\"odinger Operators.} Advanced Studies in Pure Math.~23,
Spectral and Scattering Theory and Applications (1994), 187--209.

\bibitem[JenNak2]{JN2} Jensen, A., Nakamura, S. {\em $L\sp p$-mapping properties of functions of Schr\"odinger
 operators and their applications to scattering theory.}  J.\ Math.\ Soc.\ Japan  47~(1995),  no.~2, 253--273.


\bibitem[JenNen]{JenNen} Jensen, A., Nenciu, G.
{\em A unified approach to resolvent expansions at thresholds.}
Rev.\ Math.\ Phys.\ 13 (2001), no.~6, 717--754.


\bibitem[KriSch]{KS} Krieger, J., Schlag, W. {\em On the focusing critical semi-linear wave
equation}, preprint 2005.

\bibitem[Ste1]{Stein1} Stein, E.\ M.\ {\em Topics in harmonic analysis related to the
Littlewood-Paley theory.} Annals of Mathematics Studies, No.~63
Princeton University Press, Princeton, N.J., University of Tokyo
Press, Tokyo 1970.

\bibitem[Ste2]{Stein2} Stein, E. {\em Harmonic analysis}, Princeton
University Press, Princeton 2004.

\bibitem[Yaj]{Y} Yajima, K. {\em The $W\sp {k,p}$-continuity of wave operators
for Schr\"odinger operators.}  J.\ Math.\ Soc.\ Japan  47  (1995),
no.~3, 551--581.


\bibitem[Yaj2]{Yaj2} Yajima, K. {\em Dispersive estimate for
Schr\"odinger equations with threshold resonance and eigenvalue},
preprint~2004, to appear in Comm.\ Math.\ Phys.


\end{thebibliography}

\end{document}